\documentclass[reqno]{amsart}
\usepackage[hmargin={25mm,25mm}]{geometry}
\usepackage[dvipdfmx]{graphicx,color}
\newtheorem{thm}{Theorem}[section]

\newtheorem{prop}[thm]{Proposition}
\newtheorem{cor}[thm]{Corollary}     

\theoremstyle{definition}

\theoremstyle{remark}

\title[Removing small wavenumber constraint]
{Removing small wavenumber constraints
\\
 in Side B of the Probe Method}

\author{Masaru \textsc{Ikehata}}
\address{
Professor Emeritus at Gunma University;
Professor Emeritus at Hiroshima University, Graduate School of Advanced Science and Engineering, 
Hiroshima University\\
Higashihiroshima, Japan
}
\email{ikehataprobe@gmail.com}


\begin{document}

\begin{abstract}
The Probe Method is an analytical reconstruction scheme for inverse obstacle problems utilizing the Dirichlet-to-Neumann map associated with 
the governing partial differential equation.
It consists of two distinct parts: Side A and Side B.
Both are based on the indicator sequence which is calculated from the Dirichlet-to-Neumann map acting on
"needle-like" specialized solution of the governing equation for the background medium, whose energy is concentrated on an arbitrary given needle inside.
In Side A, the limit of the indicator sequence-referred to as the indicator function-is computed before the needles touch the obstacle, and the boundary is identified as the point where this function first blows up. In contrast, Side B states the blow-up of the indicator sequence after the needles have come into contact with the obstacle.
For the Helmholtz equation, the validity of Side B has long required a small wavenumber constraint. 
This paper finally removes this long-standing restriction, establishing the method's applicability for broader cases.

\end{abstract}

\maketitle

\vskip.5cm
\noindent

\vskip.2cm
\section{Introduction}

This paper is concerned with the asymptotic behaviour of special sequences of solutions of the Helmholtz equation $\Delta v+k^2v=0$ for a fixed wavenumber $k\ge 0$.

Let $\Omega$ be a bounded domain of $\Bbb R^3$ with Lipschitz boundary \cite{Gr}.
Given $x\in\Omega$, we say a piecewise linear, non self-intersecting  curve $\sigma:[0,\,1]\longmapsto \overline{\Omega}$ is a {\it needle} with tip at $x$ if
$\sigma(0)\in\partial\Omega$, $\sigma(1)=x$ and $\sigma(t)\in\Omega$ for all $t\in\,]0,\,1[$.  
The set of all needles with tip at $x$ is denoted by $N_x$.  Hereafter, for simplicity we write $\sigma([0,\,1])=\sigma$.

Given $x\in\Omega$ and $\sigma\in N_x$ we consider a set of special sequences of solutions of the Helmholtz equation $\Delta v+k^2v=0$ in $\Omega$.
A sequence $\{v_n\}$ in $H^1(\Omega)$ is called a {\it needle sequence} for $(x,\sigma)$ if 
each term $v_n$ satisfies the Helmholtz equation $\Delta v+k^2v=0$ in $\Omega$ in the weak sense and, for any compact set $K$ of $\Bbb R^3$
with $K\subset\Omega\setminus\sigma$ it holds that
$$\displaystyle
\lim_{n\rightarrow\infty}\Vert v_n-G_k(\,\cdot\,,x)\Vert_{L^2(K)}+\Vert \nabla (v_n-G_k(\,\cdot\,,x))\Vert_{L^2(K)}=0,
$$
where $G_k=G_k(z,x)$ takes the form
$$\begin{array}{ll}
\displaystyle
G_k(z,x)=\frac{e^{ik\vert z-x\vert}}{4\pi\vert z-x\vert},
&
z\in\Omega\setminus\{x\},\,k\ge 0.
\end{array}
$$
For simplicity of description of statements, we write ${\mathcal N}(x,\sigma)$ the set of all needle sequences for $(x,\sigma)$.

It should be noted that, for any finite open cone $V$ with vertex at $x$, we have
$$\displaystyle
\int_{V\cap\Omega}\left\vert\nabla 
\text{Re}\,\left(\frac{e^{ik\vert z-x\vert}}{\vert z-x\vert}\right)\right\vert^2\,dz=\infty.
$$
Based on this property, in \cite{INew}, Lemmas 2.1 and 2.2,
it has been clarified that the needle sequence $\{v_n\}$ for $(x,\sigma)$ blows up on $\sigma$ as described below.

\begin{prop}
Given $x\in\Omega$ and $\sigma\in N_x$ let $\{v_n\}\in {\mathcal N}(x,\sigma)$.

\noindent
{\rm (a)}   Let $V$ be an arbitrary finite and open cone with vertex at $x$.  We have
$$\displaystyle
\lim_{n\rightarrow\infty}\Vert\nabla v_n\Vert_{L^2(V\cap\Omega)}=\infty.
$$

\noindent
{\rm (b)}  Given a point $z\in\sigma\setminus\{\sigma(0)\}$ let $B$ be an arbitrary open ball centered at $z$.
We have
$$\displaystyle
\lim_{n\rightarrow\infty}\Vert\nabla v_n\Vert_{L^2(B\cap\Omega)}=\infty.
$$


\end{prop}

$\quad$

{\bf\noindent Remark 1.1.}
From the fact that $G_k(\,\cdot\,,x)\in H^1_{\text{loc}}(\Omega\setminus\{x\})$ and Proposition 1.1,
for any $\{v_n\}\in {\mathcal N}(x,\sigma)$ we have the formula
$$\displaystyle
\sigma\setminus\{\sigma(0)\}=\left\{z\in\Omega\,\vert\,\lim_{n\rightarrow\infty}\Vert\nabla v_n\Vert_{L^2(B\cap\Omega)}=\infty\,\text{for any open ball $B$ centered at $z$}\right\}.
$$
The $\{v_n\}$ is generated by $\sigma$ and $\sigma$ itself is regenerated by $\{v_n\}$ via the formula above.  This is the meaning of the name `needle sequence'.

$\quad$

\noindent
From Proposition 1.1, using exactly the same argument as \cite{INew}, on p. 421, lines 10 up to 18 up, one can easily deuce

\begin{prop}
 Let $D$ be an open subset of $\Bbb R^3$ with Lipschitz boundary
and satisfy $\overline{D}\subset\Omega$.
Let $x\in\Omega$ and $\sigma\in N_x$ satisfy one of two cases {\rm (i)} and {\rm (ii)} listed below:

\noindent
{\rm (i)} $x\in\overline{D}$;

\noindent
{\rm(ii)} $x\in\Omega\setminus\overline{D}$ and $\sigma\cap D\not=\emptyset$.

\noindent
Then, for any $\{v_n\}\in {\mathcal N}(x,\sigma)$ we have
$$\displaystyle
\lim_{n\rightarrow\infty}\Vert \nabla v_n\Vert_{L^2(D)}=\infty.
$$

\end{prop}

This property was essential for the Side B of the Probe Method (e.g., \cite{IPS}, Section 4 or Section 3 in this paper) for the inverse obstacle problems governed by the Laplace and
Helmholtz equations.
However, for the Helmholtz equation (i.e., $k>0$), 
this property alone was insufficient to fully establish Side B.  The main reason is that, in Proposition 1.2
the asymptotic behabviour of the $L^2$-norm of $v_n$ over $D$ relative to that of $\nabla v_n$ was not yet well understood.
For example, in studying the inverse crack problem governed by the Laplace equation via the Probe Method, Lemma 4.3 in \cite{ICrack} established that
$$\displaystyle
\limsup_{n\rightarrow\infty}\frac{\Vert v_n\Vert_{L^2(D)}}{\Vert\nabla v_n\Vert_{L^2(D)}}<\infty,
$$
provided
$$\displaystyle
\lim_{n\rightarrow\infty}\Vert\nabla v_n\Vert_{L^2(D)}=\infty.
$$
This result holds for the $k=0$. Although the proof therein-based on the Poincar\'e inequality-also covers the case $k>0$,
this property was still insufficient to establish
Side B for the Helmholtz equation without a constraint on the size of wavenumber $k$.
Consequently, since the publication of Theorem 2.1 in \cite{INew}, the smallness constraint on $k$
 has remained an open issue; see also Theorem B in \cite{IHokkaido}.

In this paper, we take a significant step forward beyond the existing knowledge regarding the behaviour of  $\Vert v_n\Vert_{L^2(D)}$ relative to $\Vert\nabla v_n\Vert_{L^2(D)}$,
as stated below.

\begin{thm}
Let $D$ be a nonempty open subset of $\Bbb R^3$ satisfying $\overline{D}\subset\Omega$.
Assume that $D$ takes the form $D=\cup_{j=1}^N\,D_j$ with $1\le N<\infty$, where each $D_j$ is a connected component of $D$ with a Lipschitz boundary
and satisfies $\overline{D_j}\cap\overline{D_l}=\emptyset$ if $j\not=l$.
Given $x\in\Omega$ and $\sigma\in N_x$ let $\{v_n\}\in {\mathcal N}(x,\sigma)$ satisfy
\begin{equation}
\displaystyle
\lim_{n\rightarrow\infty}\,\Vert\nabla v_n\Vert_{L^2(D)}=\infty.
\tag {1.1}
\end{equation}
Then, we have
\begin{equation}
\displaystyle
\lim_{n\rightarrow\infty}
\frac{\Vert v_n\Vert_{L^2(D)}}{\Vert\nabla v_n\Vert_{L^2(D)}}=0.
\tag {1.2}
\end{equation}

\end{thm}

\noindent
Theorem 1.3 gives an answer to a question in \cite{IReview}, on page 509, Question, which is a key to an open problem goes back to \cite{INew} (2006), and
see also \cite{IRIMS}, Open problem 1.1 on page 6.

A brief outline of this paper is as follows.
Theorem 1.3 is proved in Section 2; the proof is, in short, an application of the Poincar\'e-Sobolev inequality.
In Section 3, we first review on the history of the study of Sides A and B of the Probe Method 
for an inverse obstacle problem with the impedance boundary condition governed by the Helmholtz equation with a fixed wavenumber.
We then present how Theorem 1.3 implies Side B of the Probe Method and its consequences on the characterization of the obstacle.
Section 4 is devoted to describing Side B for  the case of sound-soft obstacles, together with a review of Side A.
In the final section, a remark on a remaining problem concerning (b) of Proposition 1.1 is given.

\section{Proof of Theorem 1.3}

Fix $j$.  First, by the Poincar\'e-Sobolev inequality, e.g., see 8.11 and 8.12 on pages 220-221 in \cite{LL}, in three dimensions
we have
\begin{equation}
\displaystyle
\left\Vert v_n- \frac{1}{\vert D_j\vert}\int_{D_j}v_n\,dz\right\Vert_{L^2(D_j)}
\le L(D_j)\Vert\nabla v_n\Vert_{L^{\frac{6}{5}}(D_j)},
\tag {2.1}
\end{equation}
where $L(D_j)$ is independent of $\{v_n\}$ and $\vert D_j\vert$ denotes the Lebesgue measure of $D_j$.

Given $\epsilon>0$ set
$$\displaystyle
\sigma_{\epsilon}=\{z\in\Bbb R^3\,\vert\,\text{dist}\,(z,\sigma)\le\epsilon\}.
$$
We have $\vert \sigma_{\epsilon}\vert>0$ and $\lim_{\epsilon\rightarrow 0}\,\vert\sigma_{\epsilon}\vert=0$.
Thus $\vert D_j\setminus\sigma_{\epsilon}\vert>0$ for sufficiently small $\epsilon$.

Using the argument in the proof of Lemma in \cite{SS} (see also the proof of Proposition 2.3 in  \cite{INew}), from (2.1)
one can deduce
$$\begin{array}{ll}
\displaystyle
\left\Vert v_n- \frac{1}{\vert D_j\setminus \sigma_{\epsilon}\vert}\int_{D_j\setminus\sigma_{\epsilon}}v_n\,dz\right\Vert_{L^2(D_j)}
&
\displaystyle
\le L(D_j,\sigma_{\epsilon})
\Vert \nabla v_n\Vert_{L^{\frac{6}{5}}(D_j)},
\end{array}
$$
where
$$\displaystyle
L(D_j,\sigma_{\epsilon})=L(D_j)\left(1+\frac{\vert D_j\vert^{\frac{1}{2}}}{\vert D_j\setminus\sigma_{\epsilon}\vert^{\frac{1}{2}}}\right).
$$
Besides we have
$$\displaystyle
\left\vert\frac{1}{\vert D_j\setminus\sigma_{\epsilon}\vert}
\int_{D_j\setminus\sigma_{\epsilon}} v_n\,dz\right\vert\le\frac{1}{\vert D_j\setminus\sigma_{\epsilon}\vert^{\frac{1}{2}}}\Vert v_n\Vert_{L^2(D_j\setminus\sigma_{\epsilon})}.
$$
Thus one gets
$$\displaystyle
\Vert v_n\Vert_{L^2(D_j)}\le 
L(D_j,\sigma_{\epsilon})\Vert\nabla v_n\Vert_{L^{\frac{6}{5}}(D_j)}+\frac{\vert D_j\vert^{\frac{1}{2}}}{\vert D_j\setminus\sigma_{\epsilon}\vert^{\frac{1}{2}}}
\Vert v_n\Vert_{L^2(D_j\setminus\sigma_{\epsilon})}.
$$
This implies
\begin{equation}
\displaystyle
\frac{\Vert v_n\Vert_{L^2(D_j)}}
{\Vert\nabla v_n\Vert_{L^2(D)}}
\le L(D_j,\sigma_{\epsilon})
\frac{\Vert\nabla v_n\Vert_{L^{\frac{6}{5}}(D_j)}}{\Vert\nabla v_n\Vert_{L^2(D)}}
+\frac{\vert D_j\vert^{\frac{1}{2}}}{\vert D_j\setminus \sigma_{\epsilon}\vert^{\frac{1}{2}}}
\frac{\Vert v_n\Vert_{L^2(D_j\setminus\sigma_{\epsilon})}}{\Vert\nabla v_n\Vert_{L^2(D)}}.
\tag {2.2}
\end{equation}
Since $v_n\rightarrow G_k(\,\cdot\,,x)$ in $L^2(D_j\setminus\sigma_{\epsilon})$ and $G_k(\,\cdot\,,x)\in L^2(D_j\setminus\sigma_{\epsilon})$, from (1.1) one has
\begin{equation}
\displaystyle
\lim_{n\rightarrow\infty}\frac{\Vert v_n\Vert_{L^2(D_j\setminus\sigma_{\epsilon})}}{\Vert\nabla v_n\Vert_{L^2(D)}}=0.
\tag {2.3}
\end{equation}
Besides, we have, for all $1\le p<2$
\begin{equation}
\displaystyle
\lim_{n\rightarrow\infty}\frac{\Vert\nabla v_n\Vert_{L^{p}(D_j)}}
{\Vert \nabla v_n\Vert_{L^2(D)}}=0.
\tag {2.4}
\end{equation}
This is proved as follows.
Since $2>p$, by the H\"older inequality, we have, for all positive number $\eta$
$$\displaystyle
\Vert\nabla v_n\Vert_{L^{p}(D_j\cap\sigma_{\eta})}^{p}\le \vert D_j\cap\sigma_{\eta}\vert^{\frac{2-p}{2}}
\Vert\nabla v_n\Vert_{L^2(D_j\cap\sigma_{\eta})}^{p}
$$
and
$$\displaystyle
\Vert\nabla v_n\Vert_{L^{p}(D_j\setminus\sigma_{\eta})}^{p}\le \vert D_j\setminus\sigma_{\eta}\vert^{\frac{2-p}{2}}
\Vert\nabla v_n\Vert_{L^2(D_j\setminus\sigma_{\eta})}^{p}.
$$
Thus
$$\begin{array}{ll}
\displaystyle
\Vert\nabla v_n\Vert_{L^{p}(D_j)}^{p}
&
\displaystyle
\le\Vert\nabla v_n\Vert_{L^{p}(D_j\cap\sigma_{\eta})}^{p}+\Vert\nabla v_n\Vert_{L^{p}(D_j\setminus\sigma_{\eta})}^{p}
\\
\\
\displaystyle
&
\displaystyle
\le \vert D_j\cap\sigma_{\eta}\vert^{\frac{2-p}{2}}\Vert \nabla v_n\Vert_{L^2(D_j\cap\sigma_{\eta})}^{p}
+\vert D_j\setminus\sigma_{\eta}\vert^{\frac{2-p}{2}} \Vert\nabla v_n\Vert_{L^2(D_j\setminus\sigma_{\eta})}^{p}.
\end{array}
$$
This yields
\begin{equation}
\begin{array}{ll}
\displaystyle
\left(\frac{\Vert\nabla v_n\Vert_{L^{p}(D_j)}}
{\Vert \nabla v_n\Vert_{L^2(D)}}\right)^{p}
&
\displaystyle
\le 
\vert D_j\cap\sigma_{\eta}\vert^{\frac{2-p}{2}}\left(\frac{\Vert\nabla v_n\Vert_{L^{2}(D_j\cap\sigma_{\eta})}}
{\Vert\nabla v_n\Vert_{L^2(D)}}\right)^{p}
+\vert D_j\setminus\sigma_{\eta}\vert^{\frac{2-p}{2}}
\left(\frac{\Vert\nabla v_n\Vert_{L^{2}(D_j\setminus\sigma_{\eta})}}
{\Vert\nabla v_n\Vert_{L^2(D)}}\right)^{p}
\\
\\
\displaystyle
&
\displaystyle
\le
\vert D_j\cap\sigma_{\eta}\vert^{\frac{2-p}{2}}
+\vert D_j\setminus\sigma_{\eta}\vert^{\frac{2-p}{2}}
\left(
\frac{\Vert\nabla v_n\Vert_{L^{2}(D_j\setminus\sigma_{\eta})}}
{\Vert\nabla v_n\Vert_{L^2(D)}}
\right)^{p}.
\end{array}
\tag {2.5}
\end{equation}
Here we have $\nabla v_n\rightarrow\nabla G_k(\,\cdot\,,x)$ in $L^2(D_j\setminus\sigma_{\eta})$ and $\nabla G_k(\,\cdot\,,x)\in L^2(D_j\setminus\sigma_{\eta})$.
Thus from (1.1) one gets
$$\displaystyle
\lim_{n\rightarrow\infty}\frac{\Vert\nabla v_n\Vert_{L^{2}(D_j\setminus\sigma_{\eta})}}
{\Vert\nabla v_n\Vert_{L^2(D)}}=0.
$$
This together with (2.5) yields
$$\displaystyle
\limsup_{n\rightarrow\infty}\left(\frac{\Vert\nabla v_n\Vert_{L^{p}(D_j)}}
{\Vert \nabla v_n\Vert_{L^2(D)}}\right)^{p}
\le \vert D_j\cap\sigma_{\eta}\vert^{\frac{2-p}{2}}.
$$
Letting $\eta\rightarrow 0$, we conclude (2.4).

Then from a combination of (2.2) and (2.3) with a sufficiently small fixed $\epsilon$, and (2.4) with $p=\frac{6}{5}$
one gets
$$\displaystyle
\lim_{n\rightarrow\infty}\frac{\Vert v_n\Vert_{L^2(D_j)}}
{\Vert\nabla v\Vert_{L^2(D)}}=0.
$$
This yields (1.2).

\noindent
$\Box$

$\quad$

{\bf\noindent Remark 2.1.}
In two dimensions, instead of (2.1), we have
$$\displaystyle
\left\Vert v_n- \frac{1}{\vert D_j\vert}\int_{D_j}v_n\,dz\right\Vert_{L^2(D_j)}
\le L(D_j)\Vert\nabla v_n\Vert_{L^{1}(D_j)},
$$
where $L(D_j)$ is a positive constant being independent of $\{v_n\}$.
Hereafter, using a similar and simpler argument, one obtains Theorem 1.3 in the case $\Omega\subset\Bbb R^2$, by replacing $G_k(\,\cdot\,,x)$ with
$$\begin{array}{ll}
\displaystyle
G_k(z,x)=\frac{i}{4}H^{(1)}_0(k\vert z-x\vert), & z\in\Omega\setminus\{x\},
\end{array}
$$
where $H^{(1)}_0$ is the Hankel function of the first kind.

$\quad$

{\bf\noindent Remark 2.2.}
As can be seen from the proof, one can replace the $G_k(\,\cdot\,,x)$ and needle sequence $\{v_n\}$ forr $(x,\sigma)$
 in Theorem 1.3 with any $G(\,\cdot\,,x)\in H^1_{\text{loc}}(\Omega\setminus\{x\})$
and any sequence $\{v_n\}$ in $H^1(\Omega)$ such that $v_n\rightarrow G(\,\cdot\,,x)$ in $H^1_{\text{loc}}(\Omega\setminus\sigma)$.

$\quad$

{\bf\noindent Remark 2.3.}
The proof of (2.4) is essentially same as that of (3.5) of Proposition 3.2 in \cite{IPS4}.
Thus the point is the derivation of estimate (2.2) by using the Poincar\'e-Sobolev inequality.

\section{Implication to Side B of the Probe Method}

$\quad$

In this section, by considering an inverse obstacle problem in a bounded domain,
we show how Theorem 1.3 implies Side B of the Probe Method and its consequences on the characterization of the obstacle.

We assume that the background medium denoted by $\Omega$ be a bounded domain of $\Bbb R^3$ with $C^{1,1}$-boundary (\cite{Gr}).
Let the unknown obstacle $D$ embedded in $\Omega$ be a nonempty open set of $\Bbb R^3$ such that $\overline{D}\subset\Omega$, $\Omega\setminus\overline{D}$ is connected and
$\partial D$ is Lipschitz.  
We denote by $\nu$ both the unit outward normal to $\partial\Omega$ and $\partial D$.
Let $\lambda\in L^{\infty}(\partial D)$.  These are the starting assumptions.
Note that the $\lambda$ can be a {\it complex valued function}.

\subsection{A review on two sides of the Probe Method }

First we recall an assumption in \cite{IR}.

$\quad$

{\bf\noindent Assumption 1(\cite{IR}, Assumption 1).}  Well-posedness of the direct problem.  That is, given $F\in L^2(\Omega\setminus\overline{D})$
there  exists a unique weak solution $p\in H^1(\Omega\setminus\overline{D})$ of 
$$\left\{
\begin{array}{ll}
\displaystyle
\Delta p+k^2p=F, & x\in\Omega\setminus\overline{D},
\\
\\
\displaystyle
\frac{\partial p}{\partial\nu}+\lambda(x)p=0, & x\in\partial D,
\\
\\
\displaystyle
p=0, & x\in\partial\Omega
\end{array}
\right.
$$
and that the unique solution satisfies
$$\displaystyle
\Vert p\Vert_{L^2(\Omega\setminus\overline{D})}\le C\Vert F\Vert_{L^2(\Omega\setminus\overline{D})},
$$
where $C$ is a positive constant independent of $F$.

$\quad$

\noindent
Note that if $\text{Im}\,\lambda$ has a positive lower bound on $\partial D$, then Assumption 1 is satisfied.  This is the typical case that the obstacle surface is dissipative.
For more information about this type of boundary conditions in the exterior scattering problem, see, e.g., p.65 in \cite{CK}.

Hereafter we always assume that Assumption 1 is satisfied.
This implies that, given $f\in H^{\frac{1}{2}}(\partial\Omega)$ there exists a unique weak solution $u\in H^1(\Omega\setminus\overline{D})$
of
\begin{equation}
\left\{
\begin{array}{ll}
\displaystyle
\Delta u+k^2u=0, 
&
x\in\Omega\setminus\overline{D},
\\
\\
\displaystyle
\frac{\partial u}{\partial\nu}+\lambda(x)u=0,
&
x\in\partial D,
\\
\\
\displaystyle
u=f,
&
\displaystyle
x\in\partial\Omega.
\end{array}
\right.
\tag {3.1}
\end{equation}
Besides, as usual the Neumann data for $u$ is well-defined as the bounded linear functional 
on the $H^{\frac{1}{2}}(\partial\Omega)$,  
denoted by
$$\displaystyle
\left<\left.\frac{\partial u}{\partial\nu}\right\vert_{\partial\Omega},\,\cdot\,\right>,
$$ 
see (1.12) in \cite{IR}.  Since $u$ is uniquely determined by $f$, the map
$$\displaystyle
f\longmapsto \left<\left.\frac{\partial u}{\partial\nu}\right\vert_{\partial\Omega},\,\cdot\,\right>,
$$
is well-defined and called the Dirichlet-to-Neumann map.  It depends on $D$, $\lambda$ nonlinearly and $f$ linearly.
Similarly to the case for $\Omega\setminus\overline{D}$, the Neumann data for any weak solution $v\in H^1(\Omega)$ of the Helmholtz equation $\Delta v+k^2v=0$ in $\Omega$
is also well-defined as the bounded linear functional on $H^{\frac{1}{2}}(\partial\Omega)$,
denoted by
$$\displaystyle
\left<\left.\frac{\partial v}{\partial\nu}\right\vert_{\partial\Omega},\,\cdot\,\right>.
$$

The problem we consider in this section is

$\quad$

\noindent
{\bf Problem.} Reconstruct $D$ from the pair $(f,\frac{\partial u}{\partial\nu}\vert_{\partial\Omega})$
for infinitely many $f$, provided $\lambda$ is {\it unknown}.

$\quad$

\noindent

Let us describe a brief history of the Probe Method applied to this type of problem.
The Probe Method is an analytical approach for detecting internal discontinuities by using observational responses to "needle-like" specialized excitations.
It is based on the notion of the indicator sequence defined as below.

$\quad$

{\bf\noindent Definition 3.1(Indicator sequence, \cite{IR}, Section 5.2).}
Given $x\in\Omega$, $\sigma\in N_x$ and $\xi=\{v_n\}\in {\mathcal N}(x,\sigma)$ define
$$\displaystyle
I(x,\sigma,\xi)_n=
\text{Re}\,\left<\left.\frac{\partial v_n}{\partial\nu}\right\vert_{\partial\Omega}-\left.\frac{\partial u_n}{\partial\nu}\right\vert_{\partial\Omega},\overline{v_n}\vert_{\partial\Omega}\right>,
$$
where $u_n$ is the weak solution of (3.1) with $f=v_n\vert_{\partial\Omega}$ and $\overline{v_n}$ denotes the complex conjugate of $v_n$.
The existence of $u_n$ is a consequence of Assumption 1.

$\quad$

The Probe Method consists of two distinct parts: Side A and Side B. 
In Side A, the limit of the indicator sequence-referred to as the indicator function-is computed before the needles touch the obstacle, and the boundary is identified as the point where this function first blows up. In contrast, Side B investigates how the indicator sequence behaves after the needles have come into contact with the obstacle.

$\quad$

{\bf\noindent Side A.}
The author in \cite{IProbeScattering}(1998) and \cite{IWave}(1999) 
applied the Probe Method developed in \cite{IProbe} (1998) to the case $\lambda\equiv 0$, that is, the sound-hard obstacle and
$\lambda=\infty$, the sound-soft one.  It is based on the blow-up property of the {\it indicator function} on the obstacle surface.
Cheng-Liu-Nakamura \cite{CLN} (2003) extended his result in \cite{IProbeScattering} and \cite{IWave} to case $\lambda\not=0$.  Note that they stated their result is valid 
just for $\lambda\in L^{\infty}(\partial D)$.
However, in \cite{IHokkaido}(2006) the author pointed out that their argument does not cover this case and gave two simpler alternative proofs under the conditions
$\lambda\in C^{1}(\partial D)$ and $\partial\Omega, \partial D\in C^2$.
Recently in \cite{IR}(2022), assuming that $\lambda\in C^{0,1}(\partial D)$ and $\partial\Omega, \partial D\in C^{1,1}$,
under Assumption 1
the author obtained the lower estimate:
\begin{equation}
\displaystyle
C_1\Vert\nabla v\Vert_{L^2(D)}^2-C_2\Vert v\Vert_{L^2(D)}^2\le 
\text{Re}\,\left<\left.\frac{\partial v}{\partial\nu}\right\vert_{\partial\Omega}-\left.\frac{\partial u}{\partial\nu}\right\vert_{\partial\Omega},\overline{v}\vert_{\partial\Omega}\right>,
\tag {3.2}
\end{equation}
where $v$ is an arbitrary solution of the Helmholtz equation in $\Omega$, $u$ solves (3.1) with $f=v\vert_{\partial\Omega}$
and $C_1, C_2$ are positive constants independent of $v$.
This estimate is a special case of the more general Theorem 1.1 in \cite{IR}, which holds for a class of elliptic equations
of variable coefficients.
Here, the governing equation of {\bf Problem} is restricted to the Helmholtz equation.
As a corollary under the same assumption as above, the author established Side A of the Probe Method for {\bf Problem},
that is Theorem 5.3 in \cite{IR} described below.

\begin{thm}

\noindent
{\rm(a)}  Given $x\in\Omega\setminus\overline{D}$ let $\sigma\in N_x$ satisfy $\sigma\cap\overline{D}=\emptyset$.
Then, for any $\xi\in {\mathcal N}(x,\sigma)$ we have
$$\displaystyle
\lim_{n\rightarrow\infty}\,I(x,\sigma,\xi)_n
=I(x),
$$
where
$$\begin{array}{ll}
\displaystyle
I(x)
&
\displaystyle
=\int_{\Omega\setminus\overline{D}}\vert\nabla w_x\vert^2\,dz-k^2\int_{\Omega\setminus\overline{D}}\vert w_x\vert^2\,dz-\int_{\partial D}\text{Re}\,\lambda(z)\,\vert w_x\vert^2\,dS(z)
\\
\\
\displaystyle
&
\displaystyle
\,\,\,
+\int_D\vert\nabla G_k(z,x)\vert^2\,dz-k^2\int_{D}\vert G_k(z,x)\vert^2\,dz+\int_{\partial D}\text{Re}\,\lambda(z)\vert G_k(z,x)\vert^2\,dS(z)
\\
\\
\displaystyle
&
\displaystyle
\,\,\,
-2\int_{\partial D}\text{Im}\,\lambda(z)\,\text{Im}\,(w_x\overline{G_k(z,x)})\,dS(z),
\end{array}
$$
and $w_x=w(z)$ solves
$$
\left\{
\begin{array}{ll}
\displaystyle
\Delta w+k^2w=0, 
& 
\displaystyle
z\in\Omega\setminus\overline{D},
\\
\\
\displaystyle
\frac{\partial w}{\partial\nu}+\lambda(z)w=-\frac{\partial}{\partial\nu}G_k(z,x)-\lambda(z)G_k(z,x), 
&
\displaystyle
z\in\partial D,
\\
\\
\displaystyle
w=0,
&
\displaystyle
z\in\partial\Omega.
\end{array}
\right.
$$

\noindent
{\rm(b)}  For each $\epsilon>0$ we have
$$\displaystyle
\sup_{x\in\Omega\setminus\overline{D},\,\text{dist}\,(x, D)>\epsilon}\,\vert I(x)\vert<\infty.
$$

\noindent
{\rm(c)}  For any point $a\in\partial D$ we have
$$\displaystyle
\lim_{x\rightarrow a}\,I(x)=\infty.
$$

\end{thm}

$\quad$

\noindent
The $I(x)$ defined as the function of independent variable $x\in\Omega\setminus\overline{D}$ is called the indicator function.
In particular, (a) gives us a calculation method of the indicator function by using the indicator sequence from $\partial\Omega$ to $\partial D$.
And it describes the convergence property of the indicator sequence for $\sigma\in N_x$ satisfying $\sigma\cap\overline{D}=\emptyset$.
Needless to say, in Theorem 3.1 there is no constraint on the smallness of wavenumber.

$\quad$

{\bf\noindent Side B.}
In \cite{INew}(2005) the author formulated another side of the Probe Method,
which is based on the blow-up property of the indicator sequence when $\sigma\cap\overline{D}\not=\emptyset$.
It was first called ``Side B'' in \cite{ICrack}(2006), in which both sides of the Probe Method
were studied in an inverse crack problem governed by the Laplace equation.
A result on the Side B of the Probe Method for the {\bf Problem} was first published in \cite{IHokkaido}(2006)
provided that, roughly speaking, both $k$ and $\lambda$ are sufficiently small, see (2.7) and (2.8) of Theorem B therein.
In the seminar report \cite{IWrapping}(2005) the author noted the need to remove this constraint as an open problem.
See also \cite{IReview}, (2.32) on page 506 for the constraint when $\lambda\equiv 0$ and a brief explanation of the role in the proof of the blow-up of the indicator sequence.
However, removing these limitations has long remained an open problem.

$\quad$

\subsection{Blow-up of the indicator sequence}

Now we mention the impact of Theorem 1.3 on Side B of the Probe Method.
It is the following result, which completely removes the small wavenumber constraint combined with the upper bound of $\lambda$ from Theorem B in \cite{IHokkaido}.

\begin{thm}
Assume that: the $D$ takes the form $D=\cup_{j=1}^N D_j$ with $1\le N<\infty$, where each $D_j$ is a connected component of $D$ with a $C^{1,1}$- boundary
and satisfies $\overline{D_l}\cap\overline{D_j}=\emptyset$ if $l\not=j$; $\lambda\in C^{0,1}(\partial\Omega)$.
Let $x\in\Omega$ and $\sigma\in N_x$ satisfy one of two conditions {\rm (i)} and {\rm (ii)} listed below:

\noindent
{\rm (i)} $x\in\overline{D}$;

\noindent
{\rm (ii)} $x\in\Omega\setminus\overline{D}$ and $\sigma\cap D\not=\emptyset$.

\noindent
Then, for any $\xi=\{v_n\}\in {\mathcal N}(x,\sigma)$ we have
$$\displaystyle
\lim_{n\rightarrow\infty}\,I(x,\sigma,\xi)_n=\infty.
$$

\end{thm}

{\it\noindent Proof.}
By Theorem 1.1 in \cite{IR}, that is (3.2), we have
$$\displaystyle
C_1\Vert\nabla v_n\Vert_{L^2(D)}^2-C_2\Vert v_n\Vert_{L^2(D)}^2\le I(x,\sigma,\xi)_n.
$$
This together with Proposition 1.2 and Theorem 1.3 gives the desired conclusion.

\noindent
$\Box$

\noindent
This is the decisive result on Side B of the Probe Method for the Helmholtz equation without {\it smallness assumption} of the wavenumber $k$, even the case $\lambda\equiv 0$.
Therefore, combining Theorem 3.1, that is, Theorem 5.1 in \cite{IR}, with Theorem 3.2 above, we have finally {\it completed} the both sides of the Probe Method.
We believe the idea of the proof of Theorem 1.3 enables us to cover also other Helmholtz type PDEs/systems.
That will be reported in a subsequent paper.

\subsection{Obstacle characterization}

Hereafter, we describe how this result affect a complete characterization of unknown obstacle $D$ in terms of the asymptotic behaviour of
the indicator sequences.

First we impose the following assumption.

$\quad$

{\bf\noindent Assumption 2.}   Existence of the needle sequence, i.e., for all $x\in\Omega$ and $\sigma\in N_x$, 
the ${\mathcal N}(x,\sigma)$ is not empty.

$\quad$

\noindent
$\bullet$  It is known that, if $k^2$ is not a Dirichlet eigenvalue of $-\Delta$ in $\Omega$, then Assumption 2 is satisfied
even under the constraint on the support of $v_n$ on $\partial\Omega$.  See \cite{IWave}, Appendix, Theorem 4 together with Appendix, A.1. Remark in \cite{INew}.
Therein it is assumed that $\partial\Omega\in C^2$, however, the result is still valid for $\partial\Omega\in C^{1,1}$ using regularity theory \cite{Gr}.

\noindent
$\bullet$  If $\sigma$ is a needle given by a single segment, then one can construct an explicit needle sequence for $(x,\sigma)$ and thus ${\mathcal N}(x,\sigma)\not=\emptyset$.
See Theorem 5.1 in \cite{ICar} (also \cite{IReview} for a review).

$\quad$

The corollary below gives us a complete characterization of the obstacle $\overline{D}$ in terms of the indicator sequence, that is, via the Side B of
the Probe Method.

\begin{cor}
Assume that: the $D$ takes the form
 $D=\cup_{j=1}^N D_j$ with $1\le N<\infty$, where each $D_j$ is a connected component of $D$ with a $C^{1,1}$- boundary
and satisfies $\overline{D_l}\cap\overline{D_j}=\emptyset$ if $l\not=j$; $\lambda\in C^{0,1}(\partial\Omega)$.
Let $x\in\Omega$.  The following statements are equivalent each other.

\noindent
{\rm (a)} $x\in\overline{D}$.

\noindent
{\rm(b)} For all $\sigma\in N_x$
$$\displaystyle
\{\xi\in {\mathcal N}(x,\sigma)\,\vert\,\lim_{n\rightarrow\infty} I(x,\sigma,\xi)_n=\infty\}={\mathcal N}(x,\sigma).
$$

\noindent
{\rm(c)} For all  $\sigma\in N_x$
$$\displaystyle
\{\xi\in {\mathcal N}(x,\sigma)\,\vert\,\lim_{n\rightarrow\infty} I(x,\sigma,\xi)_n=\infty\}\not=\emptyset.
$$

\end{cor}

{\it\noindent Proof.}

\noindent
{\rm (a)} $\Rightarrow$ {\rm (b)}.  Apply Theorem 3.2.

\noindent
{\rm (b)} $\Rightarrow$ {\rm (c)}.  By Assumption 2, ${\mathcal N}(x,\sigma)\not=\emptyset$.  

\noindent
{\rm (c)} $\Rightarrow$ {\rm (a)}. This is equivalent to the statement: if $x\in\Omega\setminus\overline{D}$, then there exits a $\sigma\in N_x$
such that $\{\xi\in {\mathcal N}(x,\sigma)\,\vert\,\lim_{n\rightarrow\infty} I(x,\sigma,\xi)_n=\infty\}=\emptyset$.  This is proved as follows.
Since $\Omega\setminus\overline{D}$ is connected, there exists a $\sigma\in N_x$ such that $\sigma\cap\overline{D}=\emptyset$.
By (a) of Theorem 3.1, for any $\xi\in {\mathcal N}(x,\sigma)$, the limit $\lim_{n\rightarrow\infty}\,I(x,\sigma,\xi)_n$ exists and $-\infty<\lim_{n\rightarrow\infty}\,I(x,\sigma,\xi)_n<\infty$.  Thus one gets $\{\xi\in {\mathcal N}(x,\sigma)\,\vert\,\lim_{n\rightarrow\infty} I(x,\sigma,\xi)_n=\infty\}=\emptyset$.

\noindent
$\Box$

Besides, the following corollary  gives us a characterization of the complement of $\overline{D}$.

\begin{cor}
Under the same assumptions as Corollary 3.3, for $x\in\Omega$ the following statements are equivalent each other.

\noindent
{\rm (a)} $x\in\Omega\setminus\overline{D}$.

\noindent
{\rm(b)} There exists a $\sigma\in N_x$ such that
$$\displaystyle
\{\xi\in {\mathcal N}(x,\sigma)\,\vert\,-\infty<\lim_{n\rightarrow\infty} I(x,\sigma,\xi)_n<\infty\}={\mathcal N}(x,\sigma).
$$

\noindent
{\rm(c)} There exists a $\sigma\in N_x$ such that
$$\displaystyle
\{\xi\in {\mathcal N}(x,\sigma)\,\vert\,-\infty<\lim_{n\rightarrow\infty} I(x,\sigma,\xi)_n<\infty\}\not=\emptyset.
$$

\end{cor}

{\it\noindent Proof.}
Statements (b) and (c) are not simply the negations of those in Corollary 3.3; therefore, we describe the proof in further detail.

\noindent
{\rm (a)} $\Rightarrow$ {\rm (b)}.  By the connectedness of $\Omega\setminus\overline{D}$, one has a $\sigma\in N_x$ such that $\sigma\cap\overline{D}=\emptyset$.
Then, by {\rm (a)} of Theorem 3.1, we have
$$\displaystyle
{\mathcal N}(x,\sigma)\subset\{\xi\in {\mathcal N}(x,\sigma)\,\vert\,-\infty<\lim_{n\rightarrow\infty} I(x,\sigma,\xi)_n<\infty\}.
$$

\noindent
{\rm (b)} $\Rightarrow$ {\rm (c)}.  By Assumption 2, ${\mathcal N}(x,\sigma)\not=\emptyset$.  

\noindent
{\rm (c)} $\Rightarrow$ {\rm (a)}. This is equivalent to the statement: if $x\in\overline{D}$, then for all $\sigma\in N_x$, it holds that
$\{\xi\in {\mathcal N}(x,\sigma)\,\vert\,-\infty<\lim_{n\rightarrow\infty} I(x,\sigma,\xi)_n<\infty\}=\emptyset$.  
However, this is the direct implication of Theorem 3.2.

\noindent
$\Box$

Corollary 3.3 characterizes $\overline{D}$ through the {\it divergence} of the indicator sequence via the formula
$$\begin{array}{ll}
\overline{D}
&
=\{
x\in\Omega\,\vert\, \{\sigma\in N_x\,\vert\, {\mathcal B}(x,\sigma)={\mathcal N}(x,\sigma)\}
=N_x\}
\\
\\
&
=\{
x\in\Omega\,\vert\, \{\sigma\in N_x\,\vert\, {\mathcal B}(x,\sigma)\not=\emptyset\}
=N_x\},
\end{array}
$$
where
$$
{\mathcal B}(x,\sigma)=
\{\xi\in {\mathcal N}(x,\sigma)\,\vert\,\lim_{n\rightarrow\infty} I(x,\sigma,\xi)_n=\infty\}.
$$
This provides a contrast to the obstacle characterization by Kirsch's Factorization Method in inverse obstacle scattering, see Corollary 2.16 in \cite{KG}.

\noindent
Besides, Corollary 3.4 characterizes its complement through {\it convergence} via the formula
$$\begin{array}{ll}
\Omega\setminus\overline{D}
&
=\{
x\in\Omega\,\vert\, \{\sigma\in N_x\,\vert\, {\mathcal C}(x,\sigma)={\mathcal N}(x,\sigma)\}
\not=\emptyset\}
\\
\\
&
=\{
x\in\Omega\,\vert\, \{\sigma\in N_x\,\vert\, {\mathcal C}(x,\sigma)\not=\emptyset\}
\not=\emptyset\},
\end{array}
$$
where
$$
{\mathcal C}(x,\sigma)=
\{\xi\in {\mathcal N}(x,\sigma)\,\vert\,-\infty<\lim_{n\rightarrow\infty} I(x,\sigma,\xi)_n<\infty\}.
$$
Note that the sets ${\mathcal B}(x,\sigma)$ and ${\mathcal C}(x,\sigma)$ can be calculated from the Dirichlet-to-Neumann map
acting on ${\mathcal N}(x,\sigma)$.
It should be also noted that Corollaries 3.3 (or 3.4) gives a uniqueness theorem that:

$\quad$

${\mathcal B}(x,\sigma)$ (or ${\mathcal C}(x,\sigma)$) for all $x\in\Omega$ and $\sigma\in N_x$ uniquely determine
$D$ itself (provided $\lambda$ is also unknown).

$\quad$

{\bf\noindent Remark 3.1.}
In \cite{IHokkaido}, Corollary C, under restrictions on $k$ and $\lambda$ mentioned before, a characterization of the complement of $\overline{D}$ is given.   That is, (a) of Corollary 3.4 is equivalent to
the statement:

\noindent
{\rm(d)} There exists a $\sigma\in N_x$ such that
$$\displaystyle
\{\xi\in {\mathcal N}(x,\sigma)\,\vert\, \sup_{n}I(x,\sigma,\xi)_n<\infty\,\}\not=\emptyset.
$$
Besides, (b) implies the statement:

\noindent
{\rm(e)} There exists a $\sigma\in N_x$ such that
$$\displaystyle
\{\xi\in {\mathcal N}(x,\sigma)\,\vert\, \sup_{n}I(x,\sigma,\xi)_n<\infty\,\}={\mathcal N}(x,\sigma).
$$
\noindent
Thus by Assumption 2, (d) and (e) are equivalent each other.
Thus all the statements  (a), (b), (c) of Corollary 3.4 and (d) and (e) are equivalent each other.

$\quad$

{\bf\noindent Remark 3.2.}
The type of statement (b) of Corollary 3.3 with $k=0$ was also described in \cite{IPS}, Corollary 5.
The point of Corollary 3.3 is: (c) implies (b).

\section{Sound-soft obstacle}

The sound-soft obstacle case is excluded in the previous section.  This is the case when $\lambda=\infty$ formally.
First we give a quick review on Side A of the Probe Method for the case of sound-soft obstacles, originally established in \cite{IWave}.

In \cite{IEnclosure0}, (4.9)  we have already the basic identity
\begin{equation}
\displaystyle
\left<\left.\frac{\partial v}{\partial\nu}\right\vert_{\partial\Omega}-\left.\frac{\partial u}{\partial\nu}\right\vert_{\partial\Omega},\overline{v}\vert_{\partial\Omega}\right>
=-\Vert\nabla v\Vert_{L^2(D)}^2+k^2\Vert v\Vert_{L^2(D)}^2
-\Vert\nabla w\Vert_{L^2(\Omega\setminus\overline{D})}^2+k^2\Vert w\Vert_{L^2(\Omega\setminus\overline{D})}^2,
\tag {4.1}
\end{equation}
where $v\in H^1(\Omega)$ satisfies the Helmholtz equation in $\Delta v+k^2v=0$ in $\Omega$, $w=u-v$ and $u\in H^1(\Omega\setminus\overline{D})$
solves
\begin{equation}
\left\{
\begin{array}{ll}
\displaystyle
\Delta u+k^2u=0, 
&
x\in\Omega\setminus\overline{D},
\\
\\
\displaystyle
u=0,
&
x\in\partial D,
\\
\\
\displaystyle
u=f,
&
\displaystyle
x\in\partial\Omega,
\end{array}
\right.
\tag {4.2}
\end{equation}
with $f=v\vert_{\partial\Omega}$.  
Note that we adopt Assumption 1, but with the impedance boundary condition 
$\frac{\partial p}{\partial\nu}+\lambda(x)p=0$ on $\partial D$
replaced by the Dirichlet boundary condition $p=0$.
Then we have
$$\displaystyle
\Vert w\Vert_{H^1(\Omega\setminus\overline{D})}^2\le C_1\Vert v\Vert_{H^1(D)}^2.
$$
Besides, by (4.12) in \cite{IEnclosure0}, we have
$$\displaystyle
\Vert w\Vert_{L^2(\Omega\setminus\overline{D})}^2\le C_2\Vert v\Vert_{L^{\frac{4}{3}}(\partial D)}^2.
$$
Applying thses to (4.1), we obtain
\begin{equation}
\displaystyle
-C_3\Vert v\Vert_{H^1(D)}^2\le \text{Re}\,
\left<\left.\frac{\partial v}{\partial\nu}\right\vert_{\partial\Omega}-\left.\frac{\partial u}{\partial\nu}\right\vert_{\partial\Omega},\overline{v}\vert_{\partial\Omega}\right>
\le
-\Vert\nabla v\Vert_{L^2(D)}^2+k^2\Vert v\Vert_{L^2(D)}^2+k^2C_2\Vert v\Vert_{L^{\frac{4}{3}}(\partial D)}^2.
\tag {4.3}
\end{equation}
Here, by Theorem 1.5.1.10 in \cite{Gr} we have, for all $\epsilon\in\,]0,\,1[$
\begin{equation}
\displaystyle
\Vert v\Vert_{L^2(\partial D)}^2\le C_4(\epsilon\Vert \nabla v\Vert_{L^2(D)}^2+\epsilon^{-1}\Vert v\Vert_{L^2(D)}^2).
\tag {4.4}
\end{equation}
and since $\frac{4}{3}<2$, one has
\begin{equation}
\displaystyle
\Vert v\Vert_{L^{\frac{4}{3}}(\partial D)}^2\le C_5\Vert v\Vert_{L^2(\partial D)}^2.
\tag {4.5}
\end{equation}
Therefore, applying (4.4) with a sufficiently small $\epsilon$ and (4.5) to (4.3), we obtain
\begin{equation}
\displaystyle
-C_3\Vert v\Vert_{H^1(D)}^2\le \text{Re}\,
\left<\left.\frac{\partial v}{\partial\nu}\right\vert_{\partial\Omega}-\left.\frac{\partial u}{\partial\nu}\right\vert_{\partial\Omega},\overline{v}\vert_{\partial\Omega}\right>
\le-C_6\Vert\nabla v\Vert_{L^2(D)}^2+C_7\Vert v\Vert_{L^2(D)}^2.
\tag {4.6}
\end{equation}
From these we immediately obtain an alternative proof of Side A of the Probe Method for the sound-soft obstacle, that is,
Theorem 1 and Corollary 1 in \cite{IWave}.
 In this context, Theorem 3.1 is adapted with the following modifications:

\noindent
(i)
Replace $u_n$ with the solution of (4.2) for $f=v_n$ in the indicator sequence $I(x,\sigma,\xi)_n$ of Definition 3.1;

\noindent
(ii)
Replace $I(x)$ in Theorem 3.1 with
$$\displaystyle
I(x)=-\Vert\nabla G_k(\,\cdot\,,x)\Vert_{L^2(D)}^2+k^2\Vert G_k(\,\cdot\,,x)\Vert_{L^2(D)}^2
-\Vert\nabla w_x\Vert_{L^2(\Omega\setminus\overline{D})}^2+k^2\Vert w_x\Vert_{L^2(\Omega\setminus\overline{D})}^2,
$$
where $w_x=w$ solves 
$$
\left\{
\begin{array}{ll}
\displaystyle
\Delta w+k^2w=0, 
& 
\displaystyle
z\in\Omega\setminus\overline{D},
\\
\\
\displaystyle
w=-G_k(z,x), 
&
\displaystyle
z\in\partial D,
\\
\\
\displaystyle
w=0,
&
\displaystyle
z\in\partial\Omega.
\end{array}
\right.
$$

\noindent
(iii)
Replace the statement (c) in Theorem 3.1 with
$$
\displaystyle
\lim_{x\rightarrow a\in\partial D}I(x)=-\infty.
$$
Note that in \cite{IWave} it is assumed that both of $\partial\Omega$ and $\partial D$ are of class $C^2$, however, the argument from \cite{IEnclosure0} employed here
is valid for $C^{1,1}$ boundaries.

The new knowledge on the sound-soft obstacle case to be added in this paper is the following statement.
\begin{thm}
Assume that: the $D$ takes the form $D=\cup_{j=1}^N D_j$ with $1\le N<\infty$, where each $D_j$ is a connected component of $D$ with a $C^{1,1}$- boundary
and satisfies $\overline{D_l}\cap\overline{D_j}=\emptyset$ if $l\not=j$.
Let $x\in\Omega$ and $\sigma\in N_x$ satisfy one of two conditions {\rm (i)} and {\rm (ii)} listed below:

\noindent
{\rm (i)} $x\in\overline{D}$;

\noindent
{\rm (ii)} $x\in\Omega\setminus\overline{D}$ and $\sigma\cap D\not=\emptyset$.

\noindent
Then, for any $\xi=\{v_n\}\in {\mathcal N}(x,\sigma)$ we have
$$\displaystyle
\lim_{n\rightarrow\infty}\,I(x,\sigma,\xi)_n=-\infty.
$$

\end{thm}

This is the direct consequence of Proposition 1.2, Theorem 1.3 and the estimate derived from (4.6), that is,
$$
\displaystyle
I(x,\sigma,\xi)_n
\le-C_6\Vert\nabla v_n\Vert_{L^2(D)}^2+C_7\Vert v_n\Vert_{L^2(D)}^2.
$$

Thus, Side B of the Probe Method covers sound-soft obstacles.  
The succeeding obstacle characterizations-Corollaries 3.3 and 3.4, and Remark 3.1-under Assumption 2 remain valid with suitable changes.
In particular, in Corollary 3.3,  $\lim_{n\rightarrow\infty}I(x,\sigma,\xi)_n=\infty\Rightarrow\lim_{n\rightarrow\infty}I(x,\sigma,\xi)_n=-\infty$;
and in Remark 3.1 $\sup_nI(x,\sigma,\xi)_n<\infty\Rightarrow\inf_n I(x,\sigma,\xi)_n>-\infty$.

$\quad$

{\bf\noindent Remark 4.1.}
It follows from (4.4) for $v=v_n$ that if (1.2) is valid, then 
$$\displaystyle
\lim_{n\rightarrow\infty}\frac{\Vert v_n\Vert_{L^{2}(\partial D)}}{\Vert\nabla v_n\Vert_{L^2(D)}}=0.
$$
This together with (4.5), Proposition 1.2,  Theorem 1.3 and (4.3)  for $v=v_n$ yields also Theorem 4.1.
Although the difference lies in where inequality (4.4) is applied, the essence remains the same.

\section{A remaining problem}

Finally, we mention a minor but persistent problem that remains open.

\noindent
On Proposition 1.1, (b):  can one replace the ball $B$ with $V$, where $V$ is an arbitrary finite and open cone with vertex at $z$?
This is also a long-standing question, and resolving it would provide a more complete picture.
If the answer is affirmative, one could also describe the blow-up of $\Vert\nabla v_n\Vert_{L^2(D)}$ in Proposition 1.2 in the grazing case;
that is, the case where $(x,\sigma)$ satisfies $x\in\Omega\setminus\overline{D}$, $\sigma\cap\partial D\not=\emptyset$ and $\sigma\cap D=\emptyset$
simultaneously.  
Then, Theorems 3.2 and 4.1, in which the condition $\sigma\cap D\not=\emptyset$ in (ii) is replaced with $\sigma\cap\overline{D}\not=\emptyset$,
would be valid.  
However, as shown in the proofs of Corollaries 3.3 and 3.4, 
this modification has no effect on the characterization of the obstacle itself. For this reason, we regard it as a minor issue.

$$\quad$$

\centerline{{\bf Acknowledgment}}

The author was partially supported by Grant-in-Aid for
Scientific Research (C)(No. 24K06812) of Japan  Society for
the Promotion of Science.

$\quad$


\begin{thebibliography}{99}







\bibitem{CK}  Colton, D. and Kress, R., {\it Inverse Acoustic and Electromagnetic Scattering Theory} 4th edn, Springer, New York, 2019.




\bibitem{CLN}  Cheng, J., Liu, J. and Nakamura, G.,  Recovery of the shape of an obstacle and the boundary impedance from the far-field pattern, 
                    J. Math. Kyoto Univ., {\bf 43}(2003), No. 1, 165-186.





\bibitem{Gr} Grisvard, P.,
          \newblock  {\it Elliptic problems in nonsmooth domains}, Pitman, Boston, 1985.




\bibitem{IProbeScattering}  Ikehata, M.,
            Reconstruction of an obstacle from the scattering amplitude at a fixed frequency, Inverse Problems, {\bf 14}(1998), No.4, 949-954.




\bibitem{IProbe} Ikehata, M., 
             Reconstruction of the shape of the inclusion by boundary measurements,
             Commun. in Partial Differential Equations, {\bf 23}(1998), No.7-8, 1459-1474.





\bibitem{IWave}  Ikehata, M.,
           Reconstruction of obstacle from boundary measurements, Wave Motion, {\bf 30}(1999), No. 3, 205-223.
           




\bibitem{IEnclosure0}  Ikehata, M.,
          Reconstruction of the support function for inclusion from boundary measurements, J. Inv. Ill-Posed Probl., {\bf 8}(2000), No. 4, 367-378.




\bibitem{INew}  Ikehata, M.,
           \newblock A new formulation of the probe method and related problems,
           Inverse Problems, {\bf 21}(2005), No. 1, 413-426.






\bibitem{IWrapping}  Ikehata, M.,
           \newblock  The probe method and its applications II, {\it Seminar notes of mathematical sciences}\,(Soga, H. et. al. eds.), 
           Vol. {\bf 8}., 9-18, 2005, Ibaraki University, Mito.           




\bibitem{ICrack}  Ikehata, M.,
          \newblock Inverse crack problem and probe method,
          Cubo A Mathematical Journal, {\bf 8}(2006), No.1, 29-40.





\bibitem{IHokkaido}  Ikehata, M., 
           \newblock Two sides of probe method and obstacle with impedance boundary condition, 
           Hokkaido Math. J., {\bf 35}(2006), No. 3, 659-681.           


   
\bibitem{ICar}  Ikehata, M.,
          \newblock Probe method and a Carleman function, Inverse Problems, {\bf 23}(2007),
          No.5, 1871-1894.

   
   

\bibitem{IRIMS} Ikehata, M., The probe and enclosure methods for inverse obstacle scattering problems.  The past and present,
{\it New Developments of Functional Equations in Mathematical Analysis}, RIMS K\^{o}ky\^{u}roku, No. 1702, 2010, pp. 1-22.  
http://hdl.handle.net/2433/170012   
   
      
   
   
   



\bibitem{IR}  Ikehata, M.,
            \newblock Revisiting the probe and enclosure methods,
            Inverse Problems, {\bf 38}(2022), No.7, 075009 (33pp).






\bibitem{IReview}  Ikehata, M.,
             Extracting  discontinuity using the probe and enclosure methods,
             J. Inv. Ill-Posed Probl., {\bf 31}(2023), No. 4, 487-575.




           
\bibitem{IPS}  Ikehata, M., Integrating the probe and singular sources methods, J. Inv. Ill-Posed Probl.,
{\bf 32}(2024), No. 6, 1249-1275.



 
   
\bibitem{IPS4}  Ikehata, M., Integrating the probe and singular sources methods: IV.  IPS function for the Schr\"odinger equation, J. Inv. Ill-Posed Prbl., {\bf 34}(2026), No. 4, 605-632.
 
       
   
    
    
\bibitem{KG} Kirsch, A. and Grinberg, N.,
                  The factorization method for inverse problems, Oxford University Press: New York, 2008.





\bibitem{LL}  Lieb, L.H.  and Loss, M., {\it Analysis}, second edition, AMS, Providence, RI(2001).




\bibitem{SS}  Stanoyevitch, A. and Stegenga, D. A., Equivalence of analytic and Sobolev Poincar\'e inequalities for planar domains, Pacfic J. Math. {\bf 178}(1997), No.2, 363-375.

\end{thebibliography}
\end{document}